\documentclass[11pt]{article}
\usepackage{amssymb,geompsfi,amscd}

\def\l{{\germ l}}

\def\H{{\mathbb H}}

\def\R{{\mathbb R}}

\def\Z{{\mathbb Z}}

\def\sL{{\cal L}}

\font\l=cmr10 at 10pt
\font\ls=cmr7
\font\lss=cmr5
\font\lsy=cmsy10
\font\lsys=cmsy7
\font\lsyss=cmsy5
\font\lmi=cmmi10
\font\lmis=cmmi7
\font\lmiss=cmmi5
\font\lex=cmex10

 at 6truept
 at 6truept

\newtheorem{theorem}{Theorem}[section]
\newtheorem{lemma}[theorem]{Lemma}
\newtheorem{proposition}[theorem]{Proposition}

\newtheorem{corollary}[theorem]{Corollary}

\newcommand\heading[1]{\smallskip\noindent{\bf
#1}}
\newcommand\qed{\nopagebreak[4]\begin{flushright}\rule{0.1in}{0.1in}
\end{flushright}\pagebreak[2]}

\newcommand{\makefig}[3]{
	\begin{figure}[htbp]
        \refstepcounter{figure}
	\label{#2}
        \begin{center}
		~#3~\\
		\medskip
                {\sf Figure \thefigure.  #1}
        \end{center}
	\medskip
	\end{figure}
}

\usepackage{amsmath}
\def \link {\mathrm{link}}
\def \HH{\mathrm{H}}
\def \GL{\mathrm{GL}}
\def \Star{\mathrm{Star}}

\title{Chord Diagrams and Coxeter Links}
\smallskip
\author{Eriko Hironaka} 
\begin{document}
\maketitle

\begin{abstract}
This paper presents a construction of fibered links $(K,\Sigma)$
out of chord diagrams $\sL$.  Let $\Gamma$ be the incidence graph
of $\sL$.  Under certain conditions on $\sL$ 
the symmetrized Seifert matrix of $(K,\Sigma)$ equals the bilinear form
of the simply-laced Coxeter system $(W,S)$ associated to $\Gamma$;
and the monodromy of $(K,\Sigma)$ equals minus the Coxeter element of $(W,S)$.
Lehmer's problem is solved for the monodromy of these Coxeter links.
\end{abstract}

\section{Introduction}

A {\it chord diagram} $\sL$ is a collection of straight arcs, 
called {\it chords} on the unit disk $D\subset \R^2$ connecting mutually 
disjoint pairs of points on the boundary of $D$.  A {\it chord system}
is a chord diagram such that the chords are ordered and oriented.
Given two distinct oriented chords $\ell_1$ and $\ell_2$ define their
linking number $\link(\ell_1,\ell_2)$ to be the linking number of
their endpoints considered as oriented 0-spheres on $S^1$.

To any chord system $\sL=\{\ell_1,\dots,\ell_n\}$ associate an
$n$-dimensional inner product space $(\R^{\sL},B)$ where $\R^{\sL}$
is the $n$ dimensional vector space with basis 
$\sL$ and inner product $B$ given by 
$$
\langle \ell_i,\ell_j\rangle = 
\left \{
\begin{array}{cl}
2 & \mbox{if $i=j$}\\
\link(\ell_i,\ell_j) &\mbox{if $i<j$}\\
\link(\ell_j,\ell_i) &\mbox{if $j<i$}\\
\end{array}
\right.
$$

A fibered link is a link $K \subset S^3$, and an oriented surface $\Sigma$
whose boundary is $K$, so that $S^3 \setminus \Sigma \simeq \Sigma \times I$.
The symmetrization of the Seifert form defines an inner product $T$ on 
the vector space $\HH_1(\Sigma;\R)$. 

\begin{theorem}\label{Existence-Theorem} 
Given a chord system $\sL$ there is a fibered link $(K,\Sigma)$
together with an isomorphism
$$
\phi: \HH_1(\Sigma;\R) \rightarrow \R^{\sL}
$$
such that $B(\phi(\alpha),\phi(\beta)) = T(\alpha,\beta)$ 
for all $\alpha,\beta \in \HH_1(\Sigma;\R)$.
\end{theorem}

We will say a chord system $\sL$ is of {\it Coxeter-type} if all the
off-diagonal entries of $B$ are non-positive.  
We call $(K,\Sigma)$ a {\it Coxeter link} if it is the fibered link 
associated to a chord system $\sL$ of Coxeter-type.
Let $(W,S)$ be the simply-laced Coxeter system associated to the
incidence graph $\Gamma$ of $\sL$.
Then 
$W$ is generated by $S=\{s_1,\dots,s_n\}$, where 
each $s_i$ is the reflection on $\R^{\sL}$ 
defined by
$$
s_i (\ell_j) = \ell_j - B(\ell_i,\ell_j)\ell_i.
$$ 
The {\it Coxeter element} associated to $(W,S)$ is 
the product $c = s_1\cdots s_n \in \mbox{GL}(\R^{\sL})$.

Our second theorem relates the monodromy of a Coxeter link to the
Coxeter element of its associated system $(W,S)$.

\begin{theorem}\label{Monodromy-Theorem}  Let $(K,\Sigma)$ be a Coxeter 
link for the Coxeter system $(W,S)$, and let $h_* : \HH_1(\Sigma;\R)
\rightarrow \HH_1(\Sigma;\R)$ be the restriction of
the monodromy $h: \Sigma \rightarrow \Sigma$ of the fibration.
Then $\phi \circ {h_*} = -c \circ \phi$.
\end{theorem}

For the convenience of the reader, we review definitions and properties
of Coxeter systems in Section~\ref{Coxeter-Section} and the monodromy
of fibered links in Section~\ref{Fibered-links-Section}.  We give some 
examples and obstructions for graphs to be incidence graphs of 
chord diagrams in Section~\ref{Chord-Section}. 
In Section~\ref{Construction-Section} we prove Theorem~\ref{Existence-Theorem}
and Theorem~\ref{Monodromy-Theorem}.  Our construction
generalizes arborescent links \cite{Conway:Knots} and slalom 
links \cite{ACampo:Slalom} which apply to the case when $\Gamma$ is a
tree.  In Section~\ref{Example-Section} we give some examples of Coxeter 
links.  Finally, in Section~\ref{Lehmer-Section}, we apply our 
results to Lehmer's problem.

Coxeter links provide an easy way to construct examples of fibered links
with known monodromy.  Properties of Coxeter groups can then be used to
describe properties of the links.  For example, since iterated torus links
have finite order monodromy, if the Coxeter group is not spherical or 
affine, then the Coxeter link cannot be an iterated torus link.  

This paper was partly written while the author was supported by the 
Max-Planck-Institut of Mathematics during the summer of 2001.

\section{Coxeter Systems}\label{Coxeter-Section}

In this section we recall some properties of simply-laced Coxeter systems.  
See also \cite{Bourbaki68} or \cite{Humphreys:Coxeter} for more complete
expositions. 

Let $\Gamma$ be a finite graph with no self-loops or multiple
edges.   Let $S = \{s_1,\dots,s_n\}$ be the set of vertices.
Assume the edges of $\Gamma$ are labeled by integers $m_{i,j} \ge 3$
We say that $\Gamma$ is simply-laced if all edges are labeled
$3$.

The {\it adjacency matrix} $A$ 
of $\Gamma$ is the matrix $A = \left [ a_{i,j} \right ]$, where 
$$
a_{i,j} = \left \{
\begin{array}{cl}
1 & \qquad{\mbox{if there is an edge between $s_i$ and $s_j$}}\\
0 & \qquad{\mbox{otherwise}}
\end{array}\right . 
$$

Given a labeled graph $\Gamma$, and an ordering on the edges
$S = \{s_1,\dots,s_n\}$, there is an associated Coxeter system 
$(W,S)$, where $W$ is the finitely presented group
$$
W = \langle\ S\ : \ (s_i s_j)^{m(i,j)}\quad i,j=1,\dots,n\ \rangle.
$$
The group $W$ is called the {\it Coxeter group} associated to $\Gamma$,
and $S$ is the set of {\it Coxeter generators}.  

Let $V$ be the $n$-dimensional vector space over $\R$ with basis
$e_1,\dots,e_n$ and inner product defined by 
$$
\langle e_i, e_j \rangle = -2 \cos \left (\frac{\pi}{m_{i,j}} \right ).
$$
Let $B = [\langle e_i,e_j \rangle]$ be the associated bilinear form.

\begin{lemma} If $\Gamma$ is simply-laced, then 
$B = 2I - A$, where $I$ is the $n\times n$ identity matrix.
\end{lemma}

The {\it Coxeter representation} of $W$ 
in $\GL(V)$ is defined by
$$
s_i(e_j) = e_j - \langle e_i,e_j \rangle e_i
=\left\{ 
\begin{array}{cl}
-e_i &\qquad\mbox{if $i=j$}\\
e_j &\qquad\mbox{if $i\neq j$ and $a_{i,j}=0$}\\
 e_i+e_j &\qquad\mbox{if $a_{i,j}=1$}\\
\end{array}
\right .
$$

A Coxeter system is called {\it spherical} if its 
Coxeter group is a finite reflection group on Euclidean
space.   It is called {\it affine} if its Coxeter group
is isomorphic to a group of affine reflections.  
It is well known (see, for example, \cite{Humphreys:Coxeter} Section 4.7,
and Theorem 6.4) that $W$ is finite, and therefore $(W,S)$ is
spherical, if and only if $B$  
is positive definite, and $(W,S)$ is affine if and only
if $B$ is positive semi-definite.

The {\it Coxeter element} $c$ of $(W,S)$ is given by
$$ 
c = s_1 \cdots s_n.  
$$
Thus $c$ depends on the choice of ordering on $S$. If
$\Gamma$ is a tree, then $c$ is determined up to conjugacy 
(\cite{Humphreys:Coxeter} Proposition 3.16) and hence its
spectrum is determined by the Coxeter system.  This is not
the case if $\Gamma$ contains circuits.

The geometry of the Coxeter system is visible in the Coxeter element
(cf. \cite{ACampo:Coxeter}.)  

\begin{theorem}\label{Howlett-Theorem}{(\bf \cite{Howlett:Coxeter} Theorem 4.1)}
Let $c$ be a Coxeter element for a Coxeter system $(\Gamma, S)$.
\begin{description}
\item{(1)} $(W,S)$ is spherical if and only if
all the eigenvalues of $c$ are roots of 
unity other than $1$.  
\item{(2)} $(W,S)$ is affine if and only if $c$ has an eigenvalue
equal to 1 and all eigenvalues $c$ have modulus one.
\end{description}
\end{theorem}

For any matrix $M = [a_{i,j}]$, let $M^+$ be the strictly upper triangular
part of $M$, that is, $M^+ = [M_{i,j}]$, where
$$
M_{i,j} = \left \{
\begin{array}{lr}
a_{i,j} &\qquad \mbox{if $i < j$}\\
0 &\qquad \mbox{if $i \ge j$}
\end{array}
\right .
$$

\begin{theorem} {\bf (\cite{Howlett:Coxeter} Theorem 2.1)}
Let $(W,S)$ be a Coxeter system with bilinear form $B$.
Then $c = -U^{-1}U^t$, where
$U = I + B^+$.  
\end{theorem}

\begin{corollary}\label{Howlett-Corollary}
If $(W,S)$ is a simply-laced Coxeter system associated to the 
graph $\Gamma$, and $A$ is its adjacency matrix, then 
$c = -U^{-1}U^t$, where 
$U = I - A_+$.
\end{corollary}

\section{Monodromy of fibered links}\label{Fibered-links-Section}

Let $K$ be an oriented link in $S^3$.  Then $K$ is {\it fibered}
with fiber $\Sigma$, if 
\begin{description}
\item{(1)} $\Sigma  \subset S^3$ is an oriented surface with boundary
equal to $K$; and
\item{(2)} there is a homeomorphism 
$$
\tau :S^3 \setminus \Sigma \rightarrow \Sigma \times I,
$$
where $I$ is the open interval $(0,1)$.
\end{description}

Let $\Sigma^+ = \Sigma \times \{0\}$ and $\Sigma^- = \Sigma  \times \{1\}$.
Then $S^3 \setminus K$ is homeomorphic to $\Sigma \times I$ with
$\Sigma^-$ glued to $\Sigma^+$ by a homeomorphism
$$
h : \Sigma \rightarrow \Sigma.
$$
Here $\Sigma^-$ and $\Sigma^+$ are identified with $\Sigma$ in the obvious
way.
The induced map
$$
h_* : \HH_1(\Sigma;\Z) \rightarrow \HH_1(\Sigma;\Z)
$$
is called the {\it monodromy} of $K$, and doesn't depend on the choice 
of trivialization $\tau$.

For any loop $\gamma$ on $\Sigma$, the inclusion of $\Sigma$
in $\Sigma \times I$ induces a map
$$
\iota: \HH_1(\Sigma;\R) \rightarrow \HH_1(\Sigma^+;\R)
$$
which we will denote by $\iota(\gamma) = \gamma^+$.

Alexander duality gives a non-degenerate pairing 
$$
\HH_1(\Sigma;\R) \times \HH_1(S^3 \setminus \Sigma;\R) \rightarrow \R
$$
by linking number in $S^3$:
$$
\left (\gamma,\gamma' \right ) \mapsto \link(\gamma,\gamma').
$$
This gives a nondegenerate bilinear form 
$$
T: \HH_1(\Sigma;\R) \times \HH_1(\Sigma;\R) \rightarrow \R
$$
defined by 
$$
T(\alpha,\beta) = \link(\alpha^+,\beta) + \link(\beta^+,\alpha).
$$

Let $\gamma_1,\dots,\gamma_n$ be a basis for $\HH_1(\Sigma;\R)$.
The {\it Seifert matrix} $M$ of $K$ with respect to $\Sigma$ 
is given by $M = [\link (\gamma_i^+,\gamma_j)]$, and  $T = M + M^t$
is the symmetrization of $M$.  

The following theorem is well-known in knot theory (see,
for example, \cite{Rolfsen76}.)

\begin{theorem}\label{Seifert-Theorem} The monodromy $h_*$ of $K$ 
written with respect to the basis $\gamma_1,\dots,\gamma_n$ equals
$M^{-1}M^t$.  
\end{theorem} 

By Theorem~\ref{Seifert-Theorem} and Corollary~\ref{Howlett-Corollary}
to find a Coxeter link associated to a simply-laced Coxeter graph
$\Gamma$ with adjacency matrix $A$ it suffices to find a fibered
link whose Seifert matrix $M$ is $I - A^+$.

\section{Admissible Graphs}\label{Chord-Section}			

A chord diagram $\sL$ is a collection of straight paths
on the unit 2-disk $D$ joining pairs of points on the boundary of $D$.
The {\it incidence graph} $\Gamma$ of a chord diagram $\sL$ is 
the graph with vertices corresponding to chords and an edge
between two vertices if and only if the chords meet in the interior
of the disk.  We will call $\sL$ a
{\it realization} of $\Gamma$.  Figure \ref{lines} gives an example.
(For ease of illustration, we will draw the chords as arcs.)

\makefig{Realization of a graph.}{lines}
{\psfig{figure=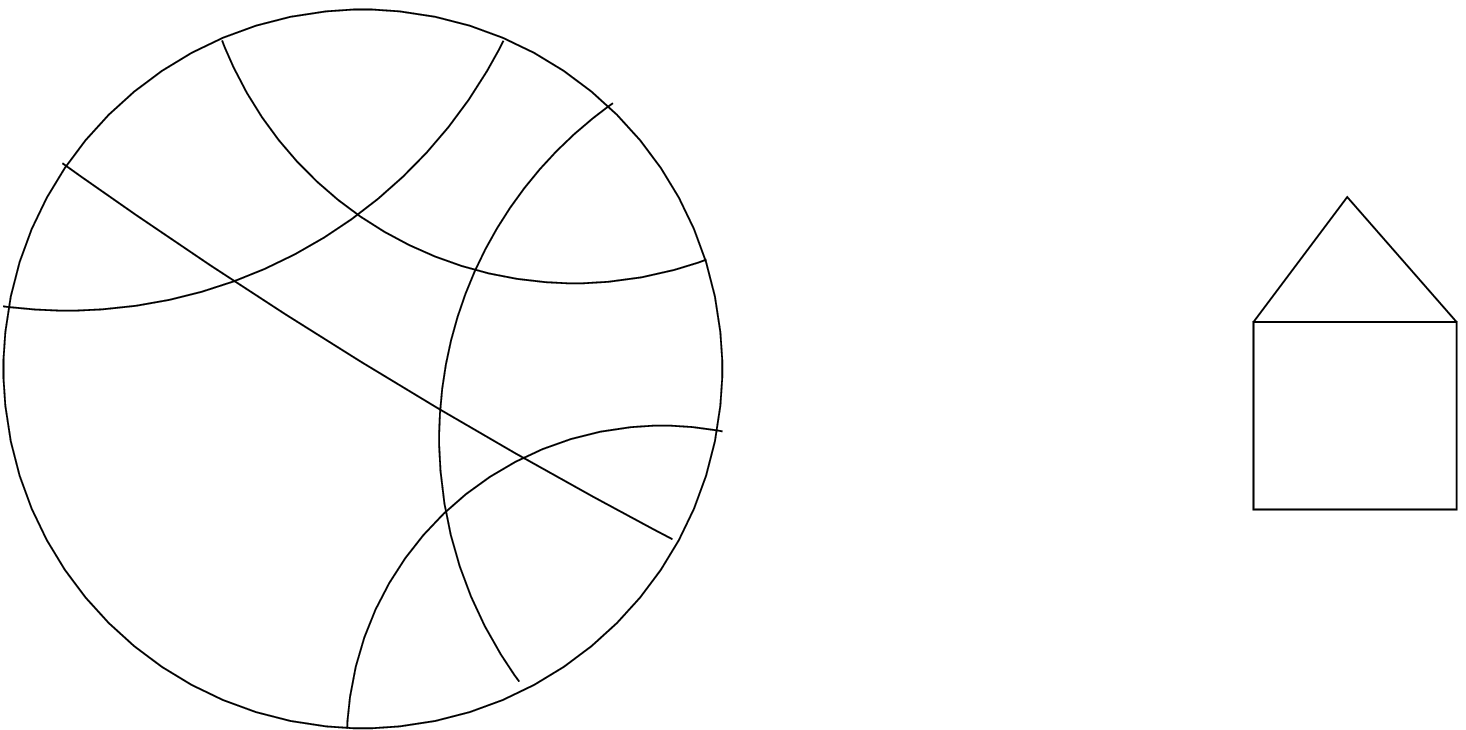,height=1.0in}}

A finite graph $\Gamma$ is {\it realizable} if it is the incidence
graph of a chord diagram.   
An ordered graph $\Gamma$ is {\it admissible} 
if there is a chord system of Coxeter-type for which 
$\Gamma$ is the incidence graph with induced ordering.
We will call two chord diagrams equivalent if they are the same up
to isotopy of the pair $(D,\sL)$. 

\begin{lemma} Given any realizable graphs 
$\Gamma_1$ and $\Gamma_2$, the join
$\Gamma_1 \vee \Gamma_2$ of the graphs at one vertex is realizable.
\end{lemma}

\heading{Proof.} A realization of any graph $\Gamma$ is equivalent
to an embedding of a union of $S^0$'s in $S^1$ one for each line in 
$\sL$.  
Let $\sL_1$ and $\sL_2$ be 
realizations of $\Gamma_1$ and $\Gamma_2$, respectively.
We can assume that
$\Gamma_1$ and  $\Gamma_2$   
correspond to a common line $\ell$ in $\sL_1$ and $\sL_2$ passing 
through the center of $D$, say horizontally as in Figure~\ref{join}. 
Furthermore, we can assume the endpoints of the arcs other than $\ell$ 
in $\sL_1$ 
lie to the left of the vertical line through the center of $D$, 
and similarly the endpoints of the arcs in $\sL_2$ other than 
$\ell$ lie to the right of the vertical line through the center of $D$.
The union of the arcs in $\sL_1$ and $\sL_2$ form a 2-embedding for 
$\Gamma_1 \vee \Gamma_2$.\nopagebreak  \qed

\makefig{Join of two realizations.}{join}
{\psfig{figure=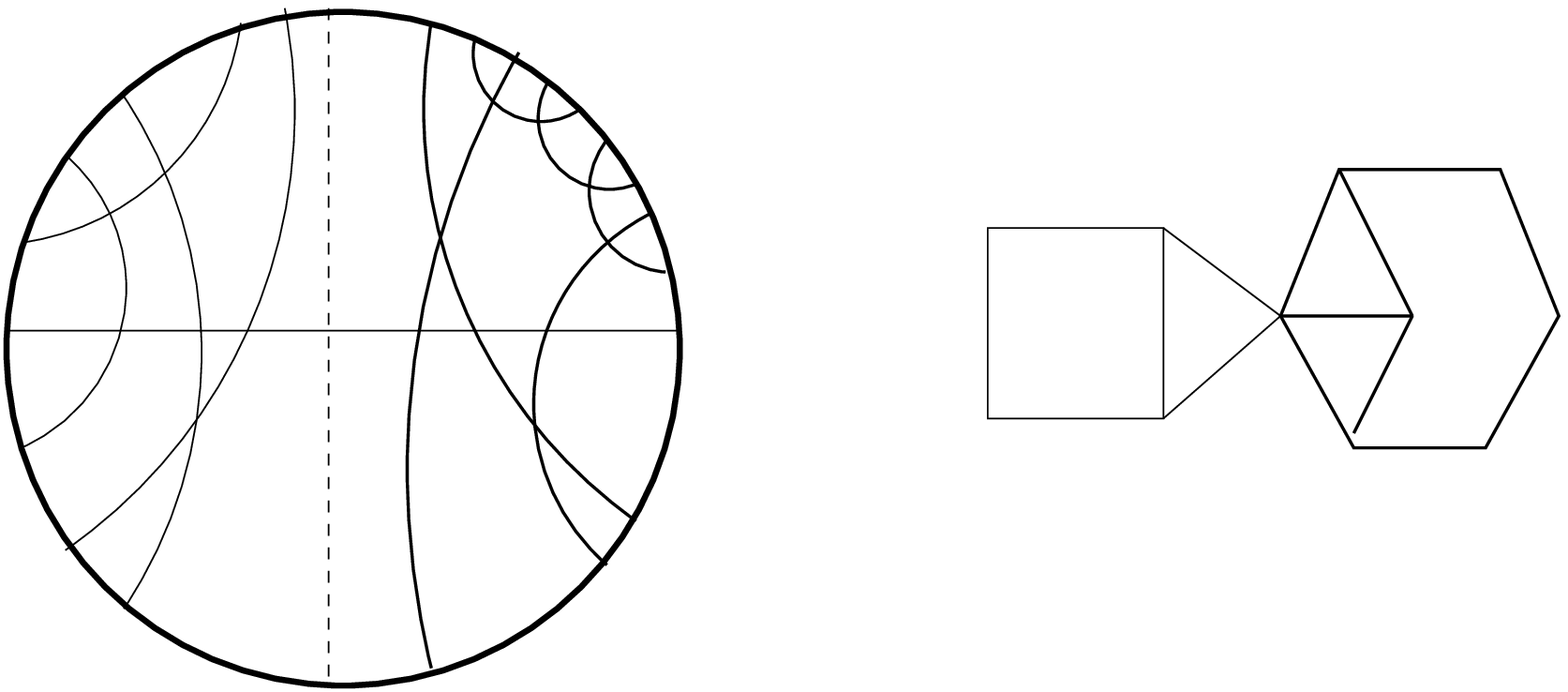,height=1.0in}}

\begin{corollary} All finite trees are realizable.
\end{corollary}

It is not hard to see that cyclic graphs, complete graphs, and 
complete bipartite graphs are realizable.  A cyclic graph has 
a realization as shown in Figure~\ref{polygon}.  

Realizations of cyclic graphs have the following property.

\begin{lemma}\label{polygon-lemma}  Up to isotopy of $(D,\sL)$ 
realizations of cyclic graphs are uniquely determined.
\end{lemma}

Let $\Gamma$ be a graph with vertices $S$.  A subgraph
$\Gamma' \subset \Gamma$ is an {\it induced subgraph} if for
some $S' \subset S$, $\Gamma'$ is the subgraph containing all
edges in $\Gamma$ whose endpoints are in $S'$.  

\makefig{5-Cycle.}{polygon}
{\psfig{figure=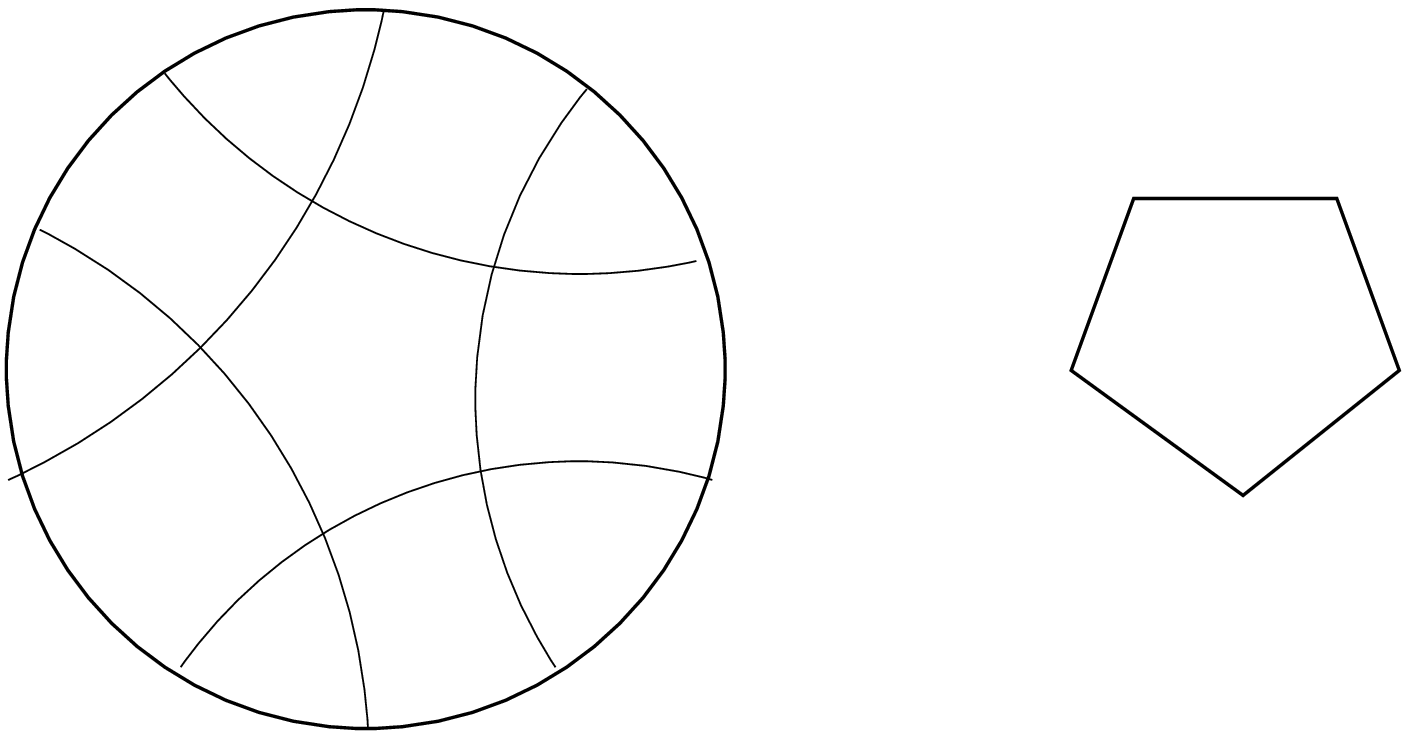,height=1.0in}}

If $S' \subset S$ is such that the induced subgraph has no edges
then we say that $S'$ is an {\it independent} set of vertices in $\Gamma$.
In order for there to exist a line in $D$ intersecting all arcs in 
an independent set $S'$, the lines in $S'$ must be parallel.  Thus, 
we have the following.

\begin{lemma} A graph 
$\Gamma$ is not realizable if there is a subset $S' \subset S$
such that
\begin{description}
\item{(1)} $S'$ contains three vertices;
\item{(2)} $S'$ is independent;
\item{(3)} there is an $s \in S$ so that for every $s'\in S'$ there
is an edge in $\Gamma$ joining $s$ and $s'$; and
\item{(4)} there is an induced cyclic subgraph in $\Gamma$ containing $S'$.
\end{description}
\end{lemma}

Figure~\ref{cube} gives an example of a non-realizable graph.

\bigskip
\makefig{Non-realizable graph.}{cube}
{\psfig{figure=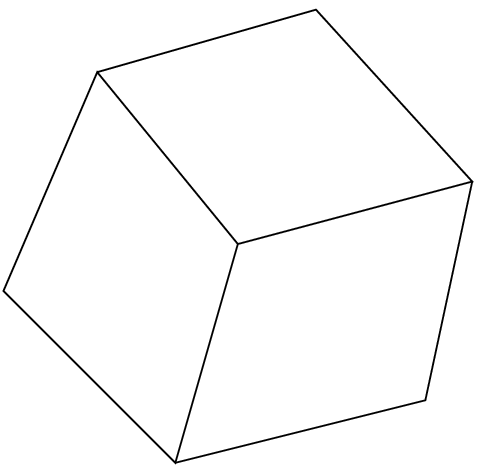,height=0.5in}}

A chord system 
$\sL = \{\ell_1,\dots,\ell_n\}$ is of
{\it Coxeter-type} if whenever $\ell_i$ and $\ell_j$
intersect in $D$, for $i < j$, the intersection looks locally as in 
Figure~\ref{ordering}.

\makefig{Orientation on a chord diagram ($i<j$).}{ordering}
{\psfig{figure=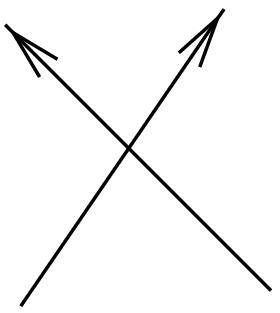,height=0.5in}}

\begin{lemma} Any chord diagram admits an ordering of Coxeter-type.
\end{lemma}

\noindent
The proof of this Lemma was communicated to me by R. Vogeler.

\heading{Proof.} 
Let $\sL$ be any realization of a chord diagram
on the unit disk in $\R^2$.  Assume that no line is horizontal.
Order the chords in $\sL$ so that for $i<j$
the slope of the line through the endpoints of $\ell_i$
has smaller slope than that of $\ell_j$.
Orient the arcs so that they are increasing
with respect to the second coordinate.
\nopagebreak
\qed

Not all orderings on a realizable graph are admissible.  
For example, given an $n$-cyclic graph, the cyclic ordering has no
Coxeter-type embedding for $n>3$.

\section{Construction}\label{Construction-Section}

Given an oriented chord diagram, we will construct an associated
fibered link.  Let 
$\sL = \{\ell_1,\dots,\ell_n\}$ be the realization of an oriented
chord diagram in the unit disk $D \subset \R^2\times \{0\} \subset \R^3$.
In $\R^3$ attach twice positively twisted bands $T_1,\dots,T_n$ to $D$ 
as in Figure~\ref{mur}, in the
order given by the ordering of the arcs, i.e., so that $T_i$ lies
over $T_j$ if $i > j$.
Let $\Sigma$ be the resulting surface, with orientation determined by
the one on $D \subset \R^3$.  Let $K_\Gamma$ be the 
oriented boundary link. 

\makefig{Murasugi sum.}{mur}
{\psfig{figure=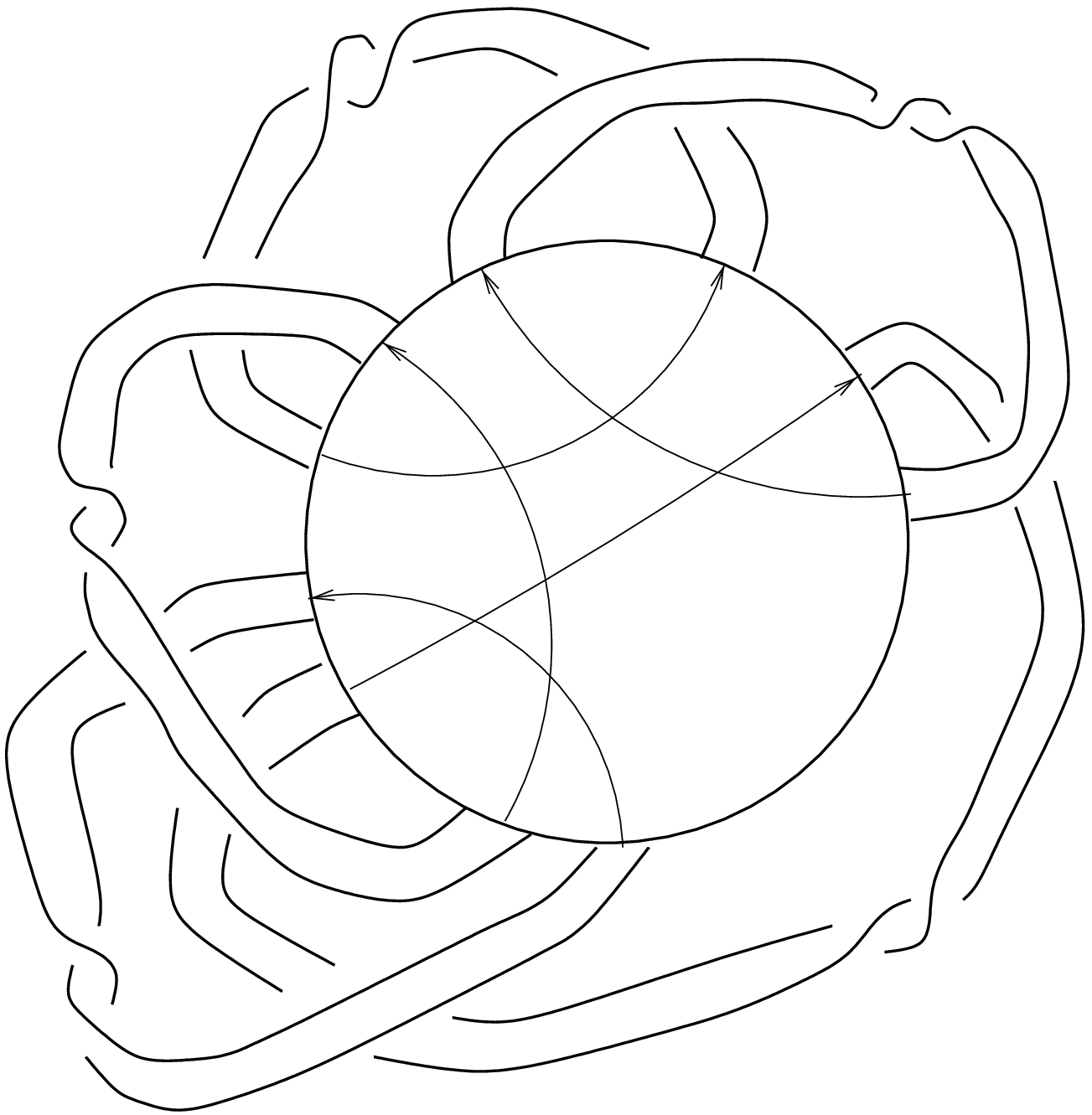,height=1.25in}}

\noindent
Then $\Sigma$ is obtained from the oriented disk $D$ by a sequence of
Murasugi sums of Hopf links.  Hence $K_\Gamma$ is a fibered link 
with fiber $\Sigma$ (\cite{Gabai:Fibered} Theorem 1). 

Extend each oriented arc $\ell_i$ to an oriented closed loop 
$\gamma_i$ going once around the corresponding attached 
handle $T_i$.  This gives a basis $\omega_1,\dots,\omega_n$
for $\HH_1(\Sigma;\R)$.  
By the construction, $\omega_i^+$ and $\omega_j$
are unlinked if $i>j$;  the positive double twist in the strand 
makes $\link(\omega_i^+,\omega_i) = 2$; and, for $i<j$, $\omega_i+$ and $\omega_j$
have linking number $\link(\ell_i, \ell_j)$.
This proves Theorem~\ref{Existence-Theorem}.

Suppose $\sL$ is of Coxeter-type.  Let $A$ be its adjacency matrix.
Then the Seifert matrix for
$(K,\Sigma)$ is given by $M = I - A^+$, and the symmetrized Seifert form
equals $2I-A$.  By Theorem~\ref{Seifert-Theorem}
the monodromy of the fibration is given by 
$$
h_* = M^{-1}M^t = U^{-1}U^t.
$$

The simply-laced Coxeter system
$(W,S)$ associated to $\Gamma$ also has bilinear form 
$$
B = 2I - A.
$$
By Theorem~\ref{Howlett-Corollary},
the Coxeter element of $(W,S)$ is given by
$$
c = -U^{-1}U^T
$$
where $U = I - A^+$.
Therefore $h_* = -c$, which proves Theorem~\ref{Monodromy-Theorem}.

We remark that our construction relied on less than the ordering of
the chord diagram.  The associated link is determined by the 
relative ordering of pairs of intersecting arcs in the chord
diagram.  We will call a chord diagram together with this information
a {\it directed chord diagram}.  Instead of an ordered incidence graph,
we obtain a directed incidence graph.

As was pointed out 
in \cite{Shi:Enumeration}, the Coxeter element of a Coxeter system
only depends up to conjugacy on the directed graph determined by 
ordered Coxeter graph.  Similarly, we can see the following from the
construction.

\begin{proposition} The fibered link associated to a chord
system only depends on the directed chord diagram.
\end{proposition}

A vertex $v$ on a directed graph is called a {\it source} (resp. {\it sink})
if all edges with one endpoint equal to $v$ point away from 
(resp. toward) $v$.  It is not hard to see that the Coxeter element
of a Coxeter system does not change its conjugacy class if a source
node is changed to a sink.  We have the following similar statement for
links constructed from chord diagrams.

\begin{proposition}\label{sink-source} The link obtained from a directed chord diagram is 
equivalent to that obtained by reordering the chord diagram so that
a source is replaced by a sink.
\end{proposition}

\heading{Proof.} Replacing a source by a sink amounts to the same as 
passing one of the twisted bands through the disk $D$
from the negative to the positive side.  Although this may change 
the isotopy type of the embedding of $\Sigma$ in $S^3 \setminus K$
it does not change the link.\qed

\section{Examples of Coxeter links.}\label{Example-Section}

This section contains some examples of Coxeter links.  

\heading{Example 1.  Trees} 

\makefig{Coxeter links and plumbing.}{tree}
{\psfig{figure=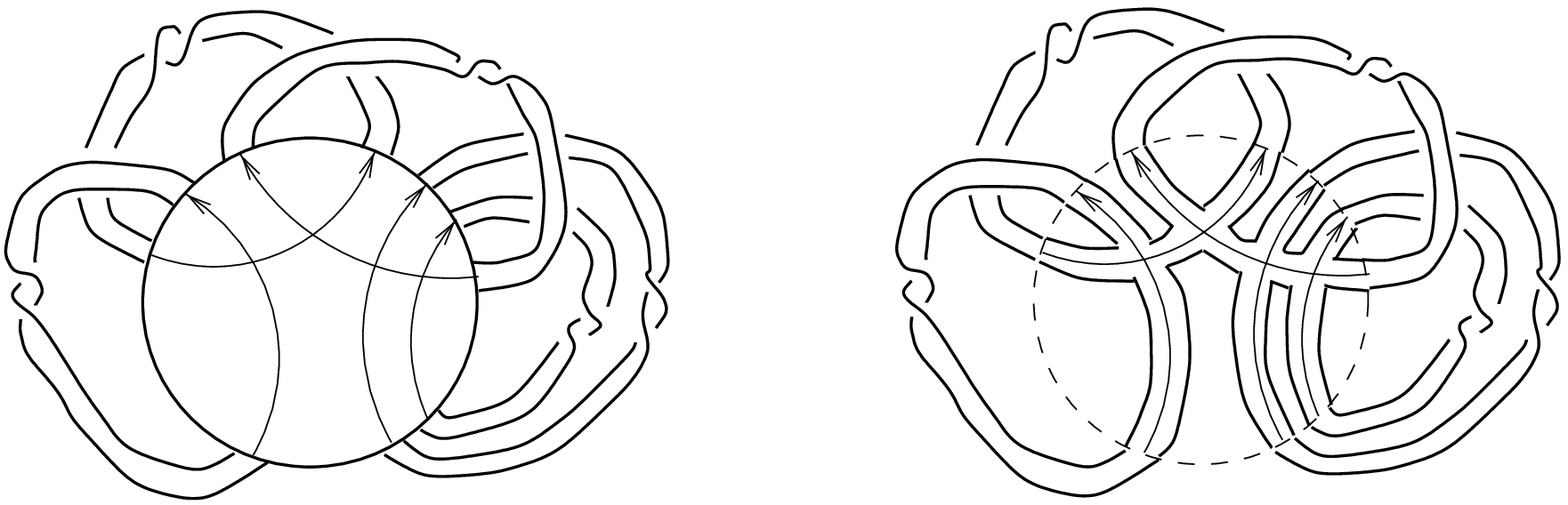,height=1in}}

When $\Gamma$ is a tree, our construction gives
arborescent links \cite{Conway:Knots}.  
This is easily seen by isotoping the 
disk to a neighborhood of the chord diagram as in
Figure~\ref{tree}.

As stated earlier in Chapter 2, the Coxeter element of a simply
laced Coxeter system $(W,S)$ doesn't
depend, up to conjugacy, on the ordering of $S$.
Visualizing the Coxeter link of a tree as a plumbing link,
one can see that the ordering of any two overlapping chords on 
a tree chord diagram can be switched by passing one of the bands
through itself.

\begin{proposition} If $\sL$ is a chord system whose incidence
graph is a tree then the associated
Coxeter link doesn't depend on the ordering on the chord
diagram.
\end{proposition}

On the other hand, there can be more than one embedding of a tree
as a positive chord system giving rise to distinct links  as
shown in Figure~\ref{treeorder}.
One sees that the link on the left has two knotted components,
while the one on the right has a component which is the unknot.
\makefig{Two embeddings of the same tree.}{treeorder}
{\psfig{figure=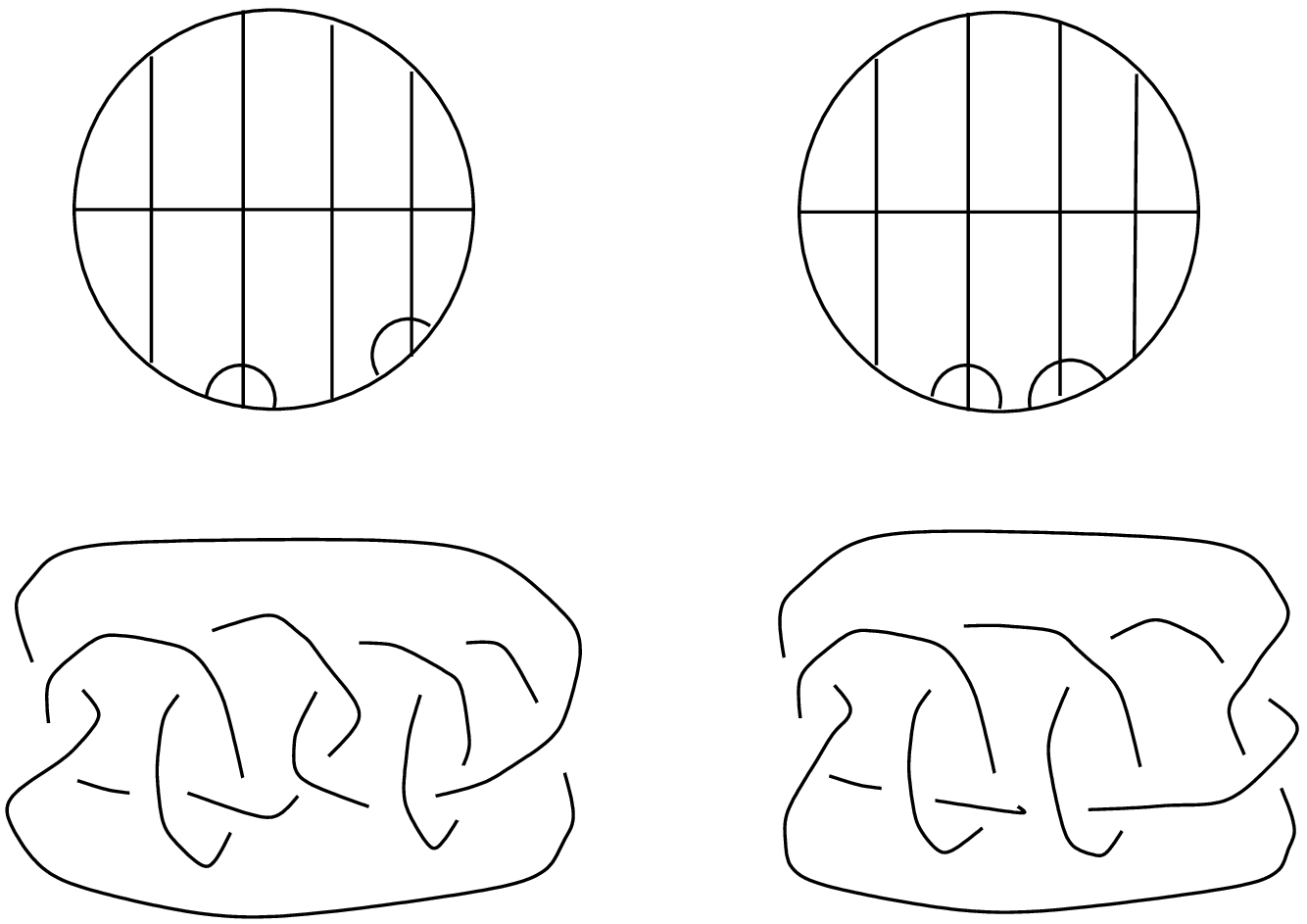,height=2in}}

\heading{Example 2. $A_n$}

The $A_n$ Coxeter graph where vertices are
numbered consecutively gives rise to the $(n+1,2)$ torus knot.

\makefig{$A_n$ gives rise to the $(n+1,2)$ torus knot.}{An}
{\psfig{figure=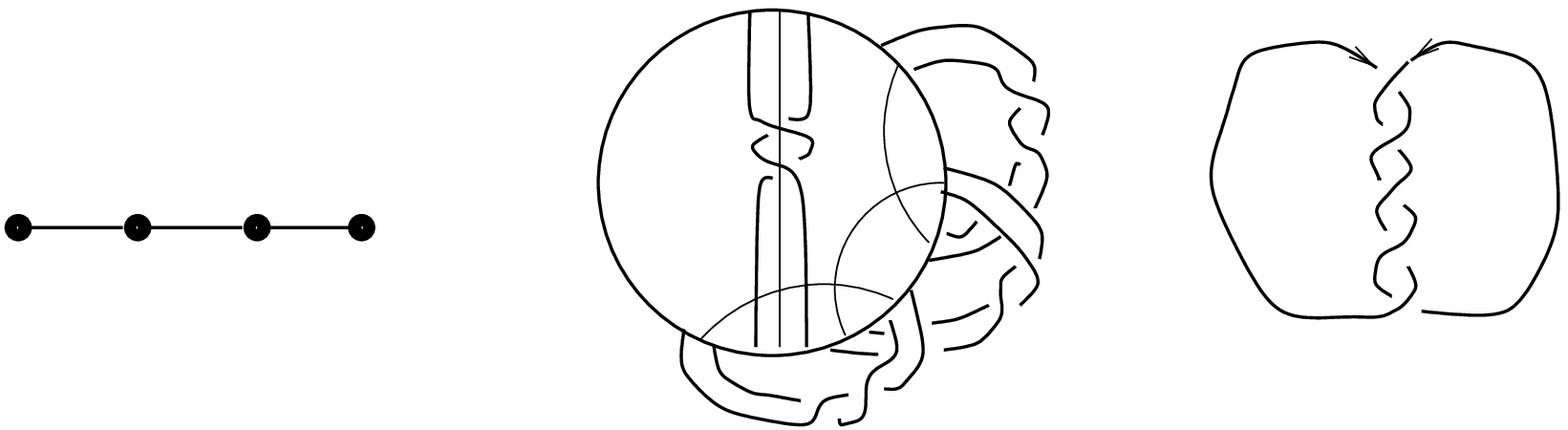,height=1.25in}}

\noindent
This can be seen inductively as follows.   A cross shaped
portion of the chord diagram, where the vertical chord
has higher index in the ordering than the horizontal one, gives rise 
to the portion of a link shown in Figure~\ref{move1}. 

\makefig{Basic transformation.}{move1}
{\psfig{figure=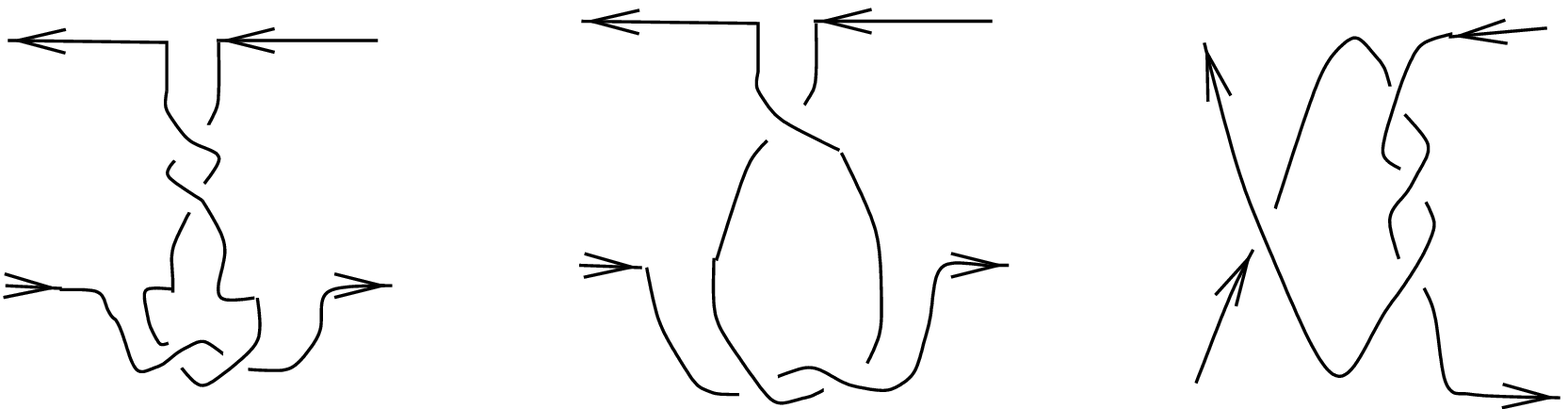,height=0.7in}}

\noindent
Thus, our claim follows by induction using the induction step
illustrated in Figure~\ref{move2}.

\makefig{Induction step.}{move2}
{\psfig{figure=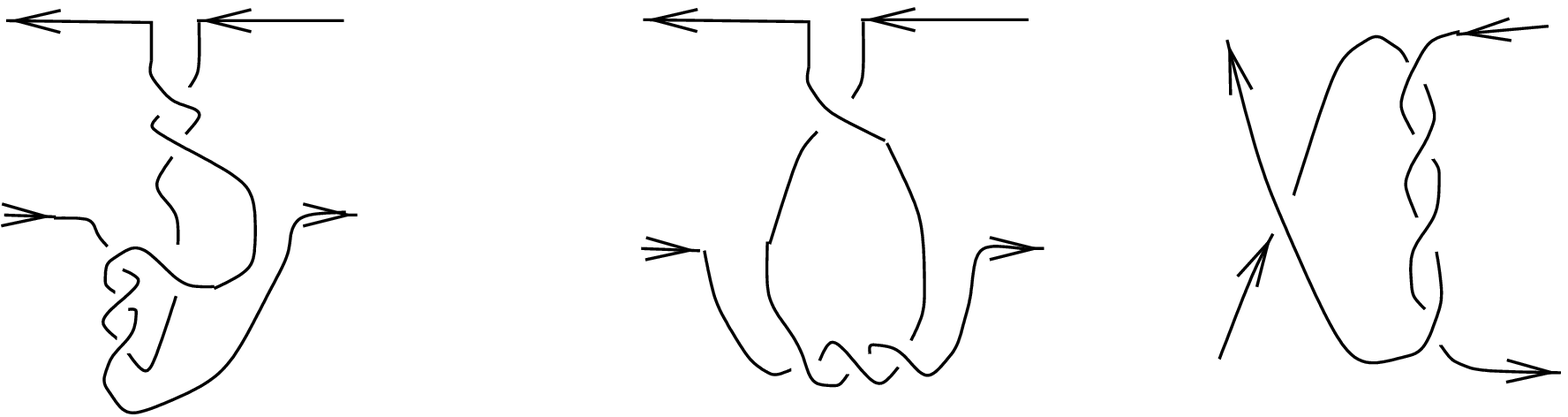,height=0.75in}}

\heading{Example 3. Star graphs}\label{Pretzel-Example} 

Let $p_1,\dots,p_k$ be positive integers.
Consider the graph $\Star(p_1,\dots,p_k)$ obtained by taking
the union of $A_{p_1},\dots,A_{p_k}$ attached at an end vertex
as in Figure~\ref{star}.  

\makefig{Star graph and its realization.}{star}
{\psfig{figure=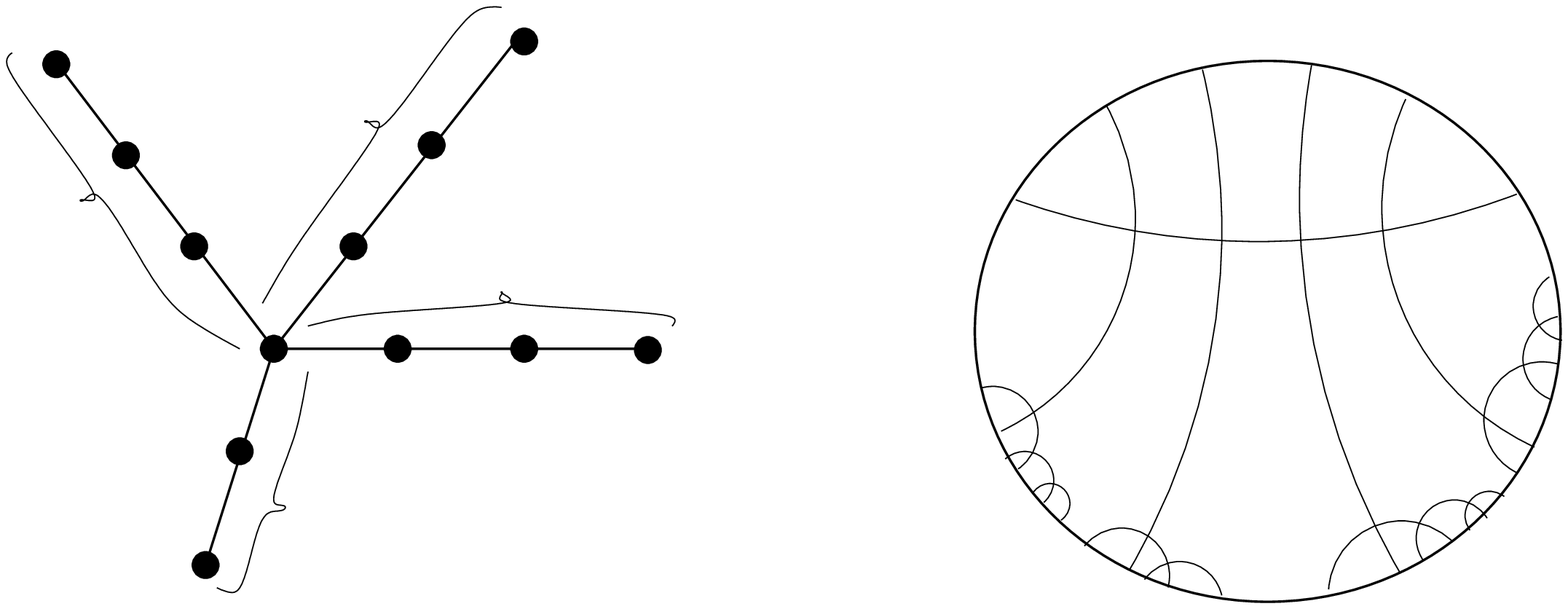,height=1in}}

\noindent 
Directing the graph so that all edges point to the multiple vertex
and using
intermediate steps shown in Figure~\ref{move1} and Figure~\ref{move2}
the reader can verify that the corresponding Coxeter link $L_{p_1,\dots,p_k}$
is a $(p_1,\dots,p_k,-1,\dots,-1)$-pretzel link, where there are $k-2$ twists
of order $-1$ (see Figure~\ref{pretzel}).  

\makefig{Coxeter link for a star graph.}{pretzel}
{\psfig{figure=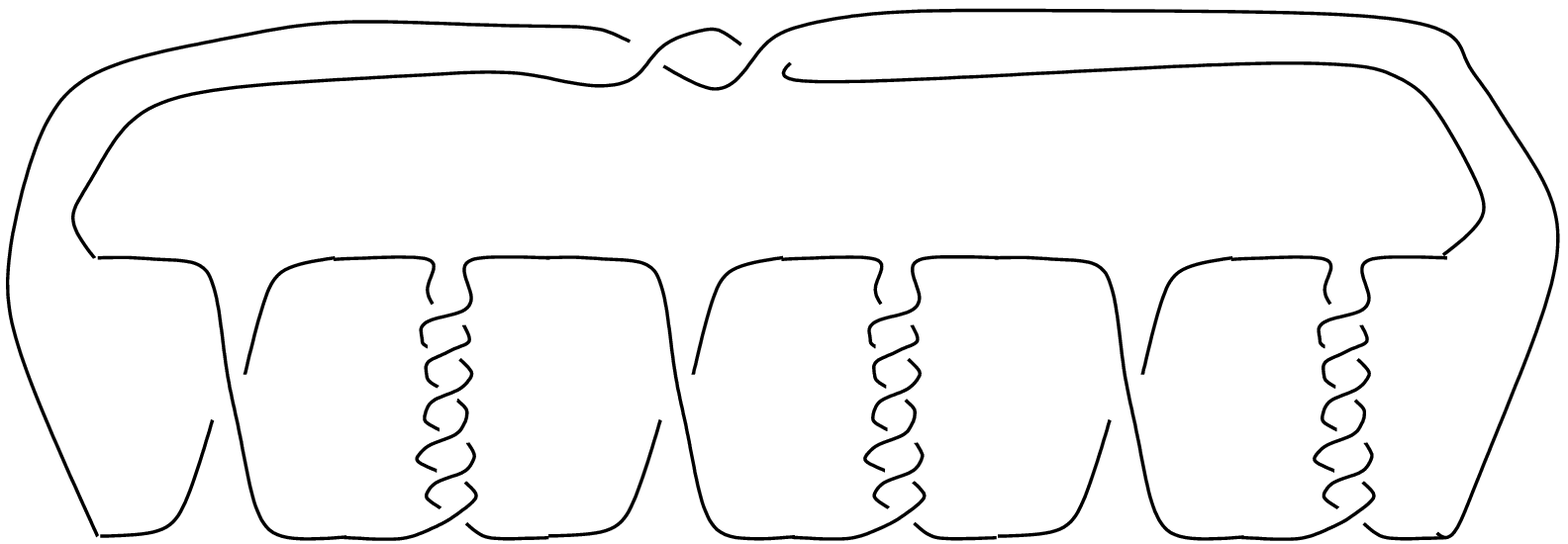,height=1in}}

Since $D_n$ is a star graph, we see that
an associated Coxeter link is the iterated torus link given 
by taking a Hopf link and replacing
one component by the $(n-1,2)$ torus link.  The 
groups $E_6$, $E_7$ and $E_8$ give rise to the $(-2,3,3)$-pretzel,
the $(-2,3,4)$-pretzel and the $(-2,3,5)$-pretzel knots, respectively.
The affine Coxeter system $E_9$ gives rise to the $(-2,3,6)$-pretzel
knot, and the hyperbolic Coxeter system $E_{10}$ gives rise to
the $(-2,3,7)$-pretzel knot.  

\heading{Example 4. $\tilde{A_n}$}\label{An-example}  

Cyclic graphs correspond to the affine Coxeter systems $\tilde{A_n}$, 
hence any Coxeter element has eigenvalue one (see Theorem 2.3), and the rest
of the eigenvalues lie on the unit circle. 

For cyclic graphs, subtleties are already exhibited for small $n$.

\makefig{Triangle systems.}{triangle}
{\psfig{figure=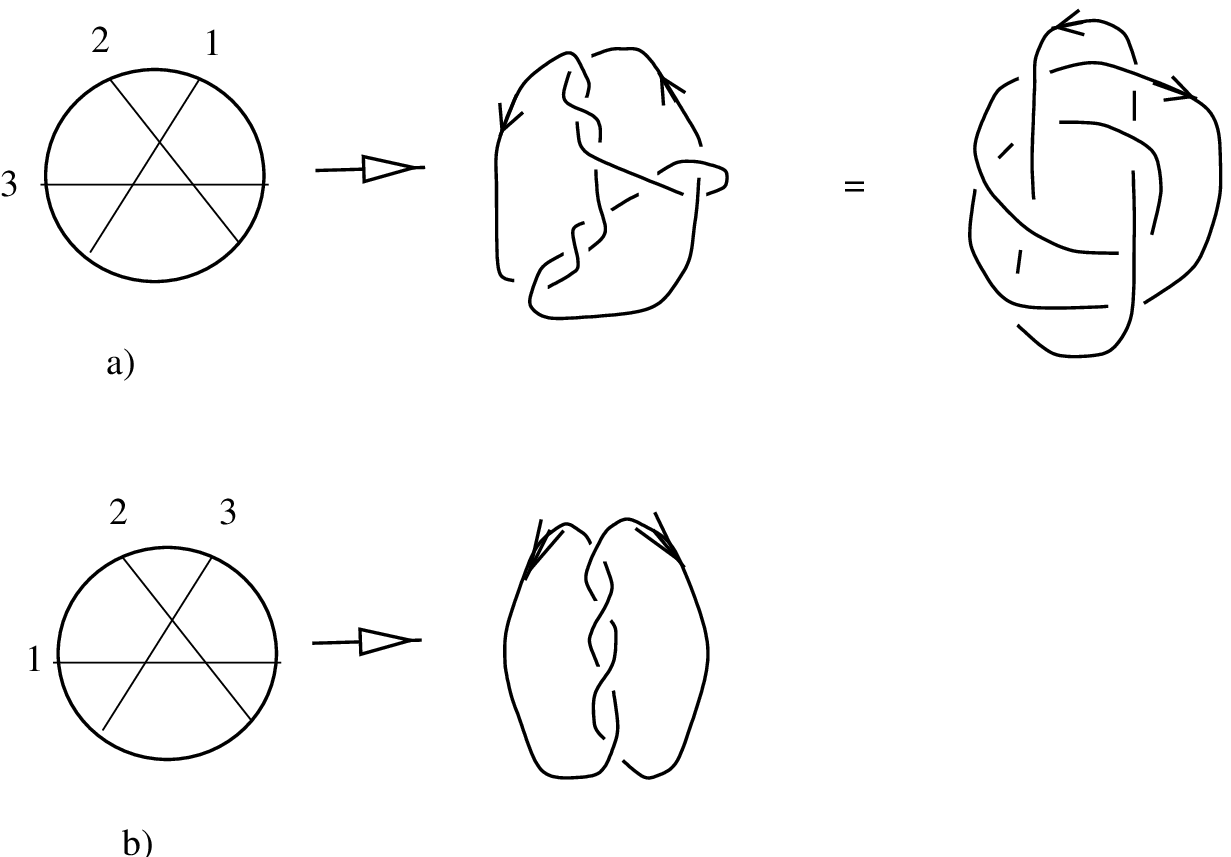,height=2.5in}}

For $n=3$, there is only one possible ordering on the
3-cycle, but there are two possible embeddings in the disk.  One of these
embeddings (Figure~\ref{triangle} a)) is of Coxeter-type and has characteristic
polynomial
$$
\Delta(t)=1 + t - t^2 - t^3 = (1+t)^2(1-t).
$$
As can be seen in the figure, this is an iterated torus link.
The other is not of Coxeter type (Figure~\ref{triangle} b)
and gives rise to the $(4,2)$ torus link $T_{4,2}$
with characteristic polynomial
$$
\Delta(t)= 1 - t + t^2 - t^3 = (1-t)(1+t^2).
$$
Note, however, that the link in Figure~\ref{triangle} b) is the Coxeter
link for the $A_3$ diagram, once we choose a different basis for 
$\H_1(\Sigma;\R)$.  

For $n=4$, there are two ordered embeddings of the cycle which are of
Coxeter-type (see Figure~\ref{square} a) and b)).  
Since the two orderings in a) and b) differ by changing a sink 
(vertex 4 in a)) 
to a source (vertex 1 in b)) the corresponding links are equivalent
by Proposition~\ref{sink-source}.  They 
equal the ($8{}^3_10$ links in Rolfsen's table \cite{Rolfsen76}).
The characteristic polynomial for both fibered links is 
given by:  
$$
\Delta(t) = 1 - 2t^2 + t^4 = (t+1)^2(t-1)^2.
$$
The embedding of the 4-cycle with cyclic ordering is not Coxeter-type.  
It gives rise to the 3-component Hopf link which has characteristic
polynomial 
$$
\Delta(t) = 1 - t - t^3 + t^4 = (1-t)^2(1+t+t^2)
$$
with respect to the fibration.
\makefig{Square systems.}{square}
{\psfig{figure=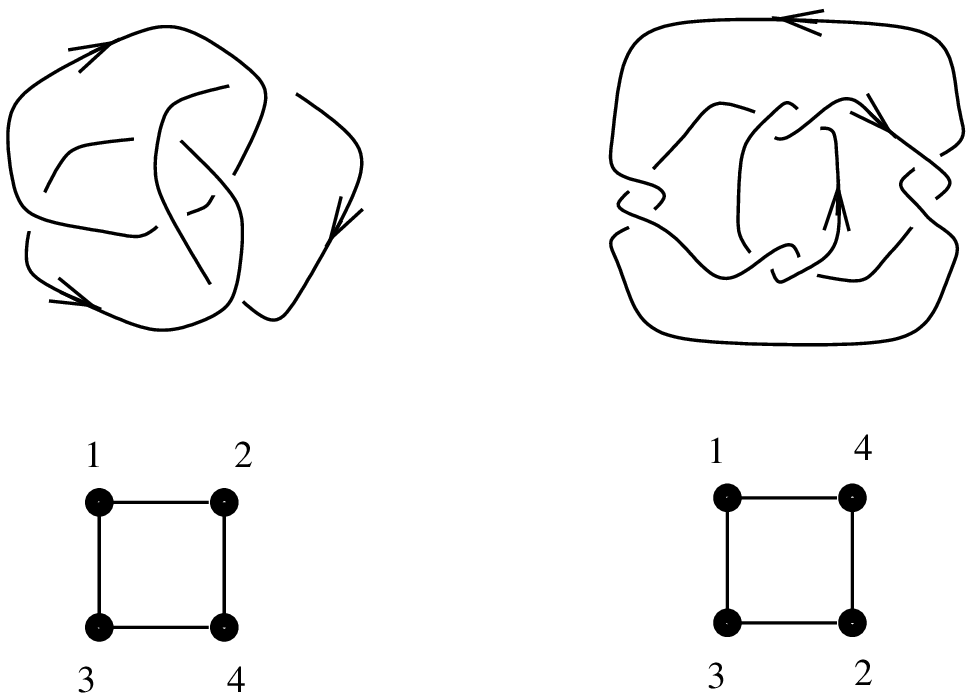,height=1.5in}}

For $n=5$, there are again two ordered embeddings of the cycle of
Coxeter-type as shown in Figure~\ref{5-cycle-graph}.
The distinct orderings give rise to the characteristic polynomials 
\begin{eqnarray*}
\Delta_1(-t) &=& 1 - t - t^4 + t^5;\quad\mbox{and}\\
\Delta_2(-t) &=& 1 - t^2 - t^3 + t^5.\\
\end{eqnarray*}
\makefig{Two orientations for the 5-cycle}{5-cycle-graph}
{\psfig{figure=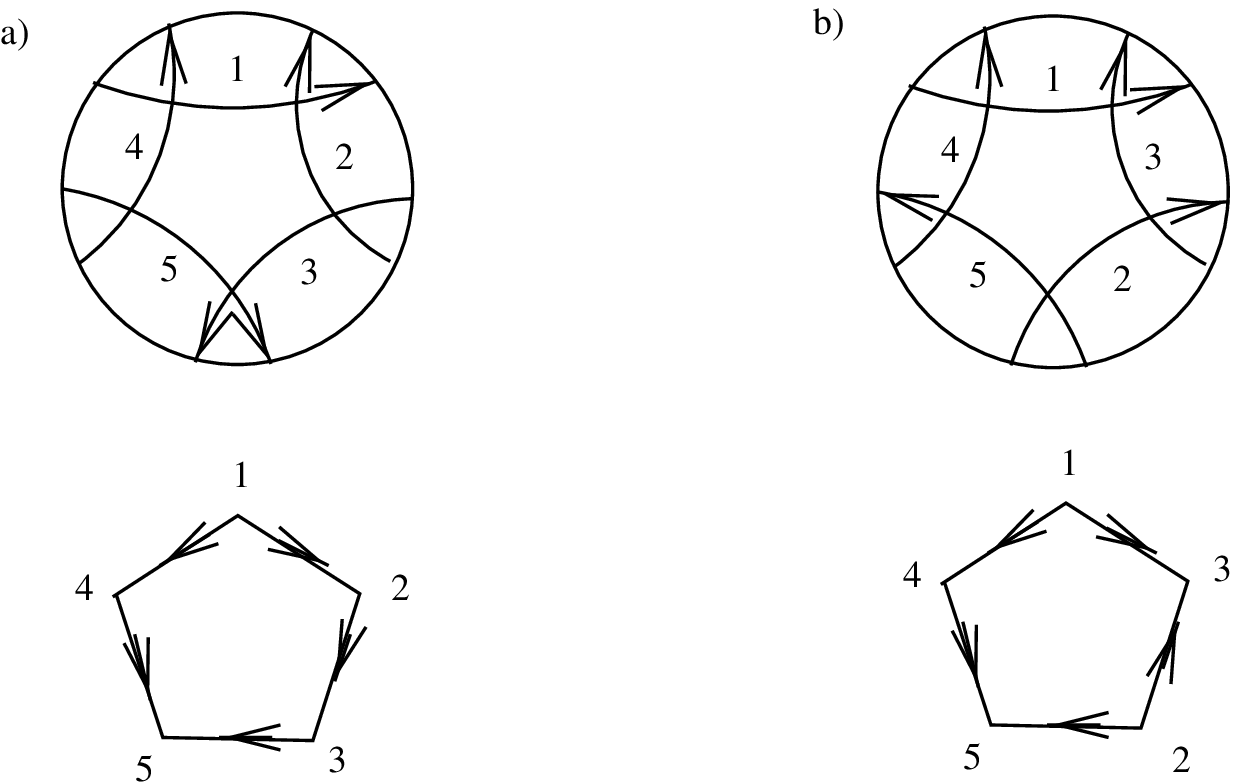,height=2in}}
The orderings give rise to the distinct links shown in Figure~\ref{5-cycle-links}.
\makefig{Two Coxeter links for the 5-cycle}{5-cycle-links}
{\psfig{figure=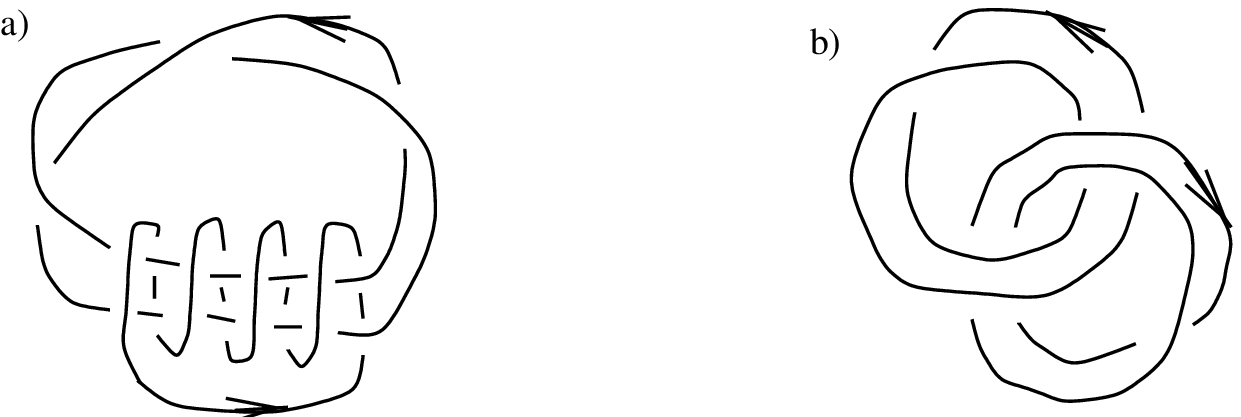,height=1in}}

\heading{5. Smallest hyperbolic Coxeter link.}

The simply-laced minimal hyperbolic Coxeter system of smallest 
dimension is a triangle with a tail, which has only one realization
as a chord diagram .    
There are three positive chord systems 
obtained by adding a chord to the positive triangle chord 
diagram in Example~\ref{An-example} in three different ways. 
By exchanging sinks and sources, however, it is possible to go from
any one of these chord systems to any other (see 
Proposition~\ref{sink-source}.) Thus, the Coxeter link
is uniquely determined and equals the knot shown in
Figure~\ref{hypert3}, which is 
the mirror of the $10{}_{145}$-knot
in Rolfsen's table \cite{Rolfsen76}.
\makefig{Smallest hyperbolic Coxeter link.}{hypert3}
{\psfig{figure=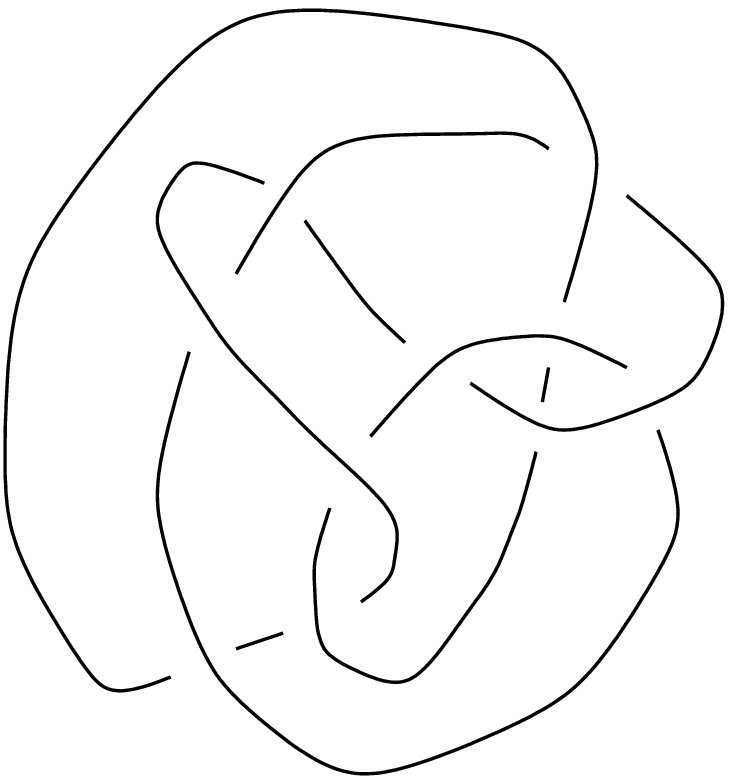,height=1in}}
The Alexander polynomial of this knot is
$$
\Delta(t) = 1 + t - 3t^2 + t^3 + t^4
$$
and its Mahler measure is: $2.36921\dots$.   

\section{Application to Lehmer's problem}\label{Lehmer-Section}

Given a polynomial $p(x) \in \Z[x]$, the {\it Mahler measure}
$\mu(p(x))$ of $p(x)$ 
is the product of the roots of $p(x)$ outside the unit 
circle.
Lehmer's problem \cite{Lehmer33} asks: For any $\delta > 0$ 
does there exist a monic integer polynomial $p(x)$ whose Mahler measure 
satisfies $1 < \mu(p(x)) < 1 + \delta$?
For degrees up to 40 (see \cite{Boyd89},\cite{Mossinghoff98}), 
the polynomial with smallest known Mahler measure 
is Lehmer's polynomial \cite{Lehmer33}
$$
p_L(x) = x^{10} + x^{9} - x^7 - x^6 - x^5 - x^4 - x^3 + x + 1.
$$

Lehmer's problem reduces to a study of monic reciprocal polynomials
\cite{Smyth70}.  
The characteristic polynomial $\Delta(x)$ of the monodromy of a fibered link is 
necessarily monic, integer, and reciprocal.  Conversely,
if $\Delta(x)$ is monic, integer, reciprocal and $\Delta(1) = \pm 1$ 
then it is the characteristic (Alexander) polynomial of fibered knot 
\cite{Burde:Alex}. For general fibered links, the restriction on 
$\Delta(1)$ does not hold.
Thus, fibered links are a natural source of examples to study Lehmer's
problem.  Kirby's Problem 5.12 in \cite{Kirby:Problems} mentions the
connection between Alexander polynomials of knots and Lehmer's
problem. Lehmer's problem is translated into a question about 
multi-variable Alexander polynomials and studied in \cite{S-W:Mahler}.  

Lehmer's polynomial appears as the Alexander polynomial of
the $(2,3,7,-1)$-pretzel knot 
\cite{Reidemeister32}.  The (-2,3,7)-pretzel knot is equivalent to the
$(2,3,7,-1)$-pretzel knot which is one of the family of pretzel
links $L_{p_1,\dots,p_k}$ (see Example~\ref{Pretzel-Example}).   
The characteristic polynomials of the monodromy $\Delta_{p_1,\dots,p_k}(x)$ 
have Mahler measure greater than or equal to that of $p_L(x)$ \cite{hironaka98}.

It is possible to answer Lehmer's problem for Coxeter links.
Let 
$$
\lambda(c) = \max\{|\alpha| \ : \ \mbox{$\alpha$ is an eigenvalue of $c$}\}.
$$
Then either $\lambda(c) = 1$ or $\lambda(c) \ge \lambda(c_0)$ 
where $c_0$ is the Coxeter element for $E_{10}$~\cite{McMullen:Coxeter}.
Thus, if $q_c$ is the characteristic polynomial of $c$,
$$
\mu(q_c(x)) \ge \mu(q_{c_0}(x))= \mu(p_L(x))
$$
for all Coxeter elements $c$.

The above discussion and the results of this paper solve Lehmer's problem
for Coxeter links.

\begin{theorem} 
If $p(x)$ be the characteristic polynomial for the monodromy of a Coxeter 
link, then $\mu(p(x)) \ge \mu(p_L(x))$.
\end{theorem} 

\bibliographystyle{math}
\bibliography{math}

\noindent
Eriko Hironaka\\
Department of Mathematics\\
Florida State University\\
Tallahassee, FL 32306\\
Email: hironaka@math.fsu.edu

\end{document}